\documentclass[a4paper]{article}

\usepackage[all]{xy}
\usepackage{amsmath, amssymb, amsthm, amsfonts}
\usepackage{graphicx}
\usepackage{longtable}
\usepackage[colorlinks=true,linkcolor=blue]{hyperref}
\usepackage{pstricks}

\newcommand{\Xfull}[1]{\overline{\mathbf{X}}_{#1}}
\newcommand{\X}[1]{\mathbf{X}_{#1}}
\newcommand{\halfspace}{H^3}
\newcommand{\Hilden}[1]{\mathbf{H}_{#1}}
\newcommand{\Braid}[1]{\mathbf{B}_{#1}}
\DeclareMathOperator{\stab}{Stab}

\newcommand{\Rorbit}{\mathcal{R}}
\newcommand{\Sorbit}{\mathcal{S}}
\newcommand{\Torbit}{\mathcal{T}}

\newcommand{\R}{\mathbb{R}}

\newcommand{\st}{\boldsymbol{\mid}}
\newcommand{\boundary}{\partial}
\newcommand{\for}{\text{for }}
\newcommand{\union}{\cup}
\newcommand{\Union}{\bigcup}
\newcommand{\intersect}{\cap}

\newtheorem{lemma}{Lemma}[section]
\newtheorem{theorem}[lemma]{Theorem}
\newtheorem{proposition}[lemma]{Proposition}
\newtheorem{corollary}[lemma]{Corollary}
\newtheorem{claim}{Claim}
\newtheorem*{claim*}{Claim}

\theoremstyle{definition}
\newtheorem{definition}[lemma]{Definition}

\def\today{\number\day \space\ifcase\month\or
  January\or February\or March\or April\or May\or June\or
  July\or August\or September\or October\or November\or December\fi
  \space\number\year}

\author{Stephen Tawn}
\title{A presentation for Hilden's subgroup of the braid group}

\begin{document}
\maketitle

\begin{abstract} 
  Consider the unit ball, $B = D \times [0,1]$, containing $n$ unknotted
  arcs $a_1, \ldots, a_n$ such that the boundary of each $a_i$ lies in 
  $D \times \{0\}$. We give a finite presentation for the mapping class
  group of $B$ fixing the arcs $\{a_1, \ldots, a_n\}$ setwise and fixing 
  $D \times \{1\}$ pointwise. This presentation is calculated using the
  action of this group on a simply-connected complex.
\end{abstract}

\section{Introduction}

Let $\halfspace$ denote the closed upper half-space of $\R^3$, let 
$a_1, a_2, \ldots, a_n \subset \halfspace$ be $n$ pairwise disjoint 
properly embedded unknotted arcs and let $a_* = a_1 \union a_2 \union 
\cdots \union a_n$.  Viewing the braid group as the mapping class group
of the punctured disc, if this disc is included in $\boundary \halfspace$ 
with $\boundary a_*$ as the punctures, one can define Hilden's group,
$\Hilden{2n}$, to be the subgroup of $\Braid{2n}$ consisting of all mapping 
classes that can be extended to $\halfspace \setminus a_*$.  Or equivalently, 
$\Hilden{2n}$ is the stabiliser of $a_*$ under the action of $\Braid{2n}$
on $0,2n$--tangles.  

Hilden\cite{Hilden} found generators for a similar
subgroup of the braid group of a sphere. For any given braid $b$
multiplying on either the left or the right by elements of $\Hilden{2n}$
preserves the plat closure, ie plat closure is constant on each double
coset.  Birman\cite{Birman} showed that if two braids have the same plat
closure then they can be related by a sequence of these double coset
moves and stabilisation moves that changes the braid index by 2.

We calculate a presentation for $\Hilden{2n}$ using the action of this
group on a cellular complex. Hatcher--Thurston\cite{HatcherThurston}, 
Wajnryb\cite{Wajnryb1,Wajnryb3,Wajnryb2}, Laudenbach\cite{Laudenbach}, 
etc used the same method to calculate presentations for mapping class
groups. We start in Section~\ref{method} by outlining this method. A
similar but more general method is given by Brown \cite{Brown}.

In Section~\ref{Xfull} we define a simply-connected complex $\Xfull{n}$.
In Section~\ref{X} we remove some of the edges and faces of this complex
resulting in a new complex which remains simply-connected but gives a
simpler presentation.  This presentation is calculated in
Section~\ref{calc} and then used to calculate a presentation with
generators similar to those found by Hilden.

\section{The method} \label{method}

Suppose that $X$ is a connected simply-connected cellular complex in which
each edge is uniquely determined by its endpoints, that $G$ is a group
acting cellularly on the right of $X$, and that this action is
transitive on the vertex set $X^0$. Pick a vertex $v_0 \in X^0$ as a
basepoint and let $H$ denote its stabiliser in $G$, ie 
$H = \{g \in G \st g \cdot v_0 = v_0 \}$.  Suppose that $H$ has a
presentation with generating set $S_0$ and relations $R_0$, ie
$H = \langle S_0 | R_0 \rangle$.

Given vertices $u,v \in X^0$ such that $\{u,v\}$ is the boundary of an
edge of $X$ we will write $(u,v)$ for this (oriented) edge.
Given a sequence $v_1, v_2, \ldots, v_k$ of vertices such that either
$v_i = v_{i+1}$ or $(v_i, v_{i+1})$ forms an edge we will write
$(v_1, v_2, \ldots, v_k)$ for the path traversing these edges.
Whenever $v_i = v_{i+1}$ we shall say that $v_i$ is a stationary point.

Let $E$ denote the set of all oriented edges starting at $v_0$, so $H$
acts on $E$.  Suppose that $\{e_\lambda\}_{\lambda\in\Lambda}$ is a
set of representatives for the $H$--orbits of the edges in $E$, ie 
$E = \Union_{\lambda \in \Lambda} H e_\lambda$ and 
$H e_\lambda = H e_{\lambda'}$ only if $\lambda = \lambda'$.  Since the
action of $G$ is transitive on $X^0$ we can find $r_\lambda \in G$
such that $e_\lambda = (v_0, v_0 \cdot r_\lambda)$.  Let
$S_1 = \{r_\lambda\}_{\lambda\in\Lambda}$.

The edges $\{e_\lambda\}_{\lambda\in\Lambda}$ also form a set of 
representatives for the edge orbits of the $G$--action on $X$.  To see
this suppose that two of these edges lie in the same $G$--orbit, ie
$(v_0,v) = (v_0, u) \cdot g$.  Then we have that $v_0 = v_0 \cdot g$
therefore $g \in H$.

Suppose that $\{f_\mu\}_{\mu\in M}$ is a set of representatives for the
$G$--orbits of the faces of $X$.  Since the action is transitive on
$X^0$, we may assume that the boundary of each face $f_\mu$ contains
the vertex $v_0$.

\begin{definition}
  An \emph{h-product of length} $k$ is a word of the form
  \[ 
    h_{k+1}\ r_{\lambda_k} h_k\ r_{\lambda_{k-1}} h_{k-1}\ 
                                             \cdots\ r_{\lambda_1} h_1
  \]
  where each $\lambda_i \in \Lambda$ and each of the $h_i$ are words in
  $H$. To each h-product we can associate an edge path 
  $P = (v_0, v_1, \ldots, v_k)$ in $X$ starting at $v_0$ then visiting
  the vertices $v_1 = v_0 \cdot r_{\lambda_1} h_1$, 
  $v_2 = v_0 \cdot r_{\lambda_2} h_2\ r_{\lambda_1} h_1$, etc. This means
  that the edge $(v_{i-1}, v_i)$ is in the orbit of
  $(v_0, v_0 \cdot r_{\lambda_i})$.  Given any edge path starting at
  $v_0$ we can choose an h-product to represent it.
\end{definition}

We can now choose the following three sets of relations.
\begin{itemize}
\item[$R_1$:] For each edge orbit representative $e_\lambda$ pick a
  generating set $T$ for the stabiliser of this edge, ie
  $\langle T \rangle = \stab_G(v_0) 
                             \intersect \stab_G(v_0 \cdot r_\lambda)$.
  For each $t \in T$ we have the relation 
  $r_{\lambda} t r_{\lambda}^{-1}=h$ for some word $h \in H$.

\item[$R_2$:] For each $e_\lambda$ we have a relation 
  $r_{\lambda'} h\ r_\lambda = h'$ where the LHS is a choice of h-product
  for the path $(v_0,v_0 \cdot r_\lambda, v_0)$ and $h'$ is some word
  in $H$.

\item[$R_3$:]  For each face orbit representative $f_\mu$ with boundary
  $(v_0, v_1, \ldots, v_{k-1}, v_0)$ choose an h-product representing
  this path and a word $h \in H$ such that 
  $r_{\lambda_k} h_k \cdots r_{\lambda_1} h_1 = h$.
\end{itemize}

\begin{theorem}\label{mainthm}
  The group $G$ has the following presentation.
  \[ G = \langle S_0 \union S_1
                     | R_0 \union R_1 \union R_2 \union R_3 \rangle \]
\end{theorem}

\begin{corollary}
  Suppose that $H$ is finitely presented, that the number of edge and 
  face orbits is finite and that each edge stabiliser is finitely generated.
  Then $G$ has a finite presentation.
\end{corollary}

We prove Theorem \ref{mainthm} in several steps.

\begin{claim} The set $S_0 \union S_1$ generates $G$. \end{claim}
\begin{proof}
  Given any $g \in G$, let $v = v_0 \cdot g$.  Now as $X$ is connected
  there is an edge path connecting $v_0$ to $v$.  Choose an h-product
  $g_1 = h_{k+1}\ r_{\lambda_k} h_k \cdots r_{\lambda_1} h_1$ representing
  this path. Then $v_0 \cdot g g_1^{-1} = v_0$ so $g = h g_1$ for some
  $h \in H$.
\end{proof}

\begin{claim} \label{claim0+1}
  If two h-products, $p_1$ and $p_2$, give rise to the same path and
  are equal in $G$ then they are equivalent modulo $R_0 \union R_1$.
\end{claim}
\begin{proof}
  Because $p_1$ and $p_2$ represent the same path they must have equal
  length.  Suppose that
  $p_1 = h_{k+1}\ r_{\lambda_k} h_k \cdots r_{\lambda_1} h_1$ and
  $p_2 = f_{k+1}\ r_{\lambda'_k} f_k \cdots r_{\lambda'_1} f_1$.
  Clearly, if the two h-products are of length 0 then they are both
  words in $H$ and so are equivalent modulo $R_0$.  Now suppose that
  $k \ne 0$.  The fact that $p_1$ and $p_2$ represent the same path
  means that
  \[ 
    (v_0,\, v_0 \cdot r_{\lambda_1} h_1,\, 
           v_0 \cdot r_{\lambda_2} h_2\ r_{\lambda_1} h_1,\, \ldots)
    = (v_0,\, v_0 \cdot r_{\lambda'_1} f_1,\, 
           v_0 \cdot r_{\lambda'_2} f_2\ r_{\lambda'_1} f_1,\, \ldots),
  \]
  therefore
  \[ (v_0, v_0 \cdot r_{\lambda_1}) 
               = (v_0, v_0 \cdot r_{\lambda'_1}) \cdot f_1 h_1^{-1}. \]
  So $\lambda_1 = \lambda'_1$ and $f_1 h_1^{-1}$ is in the stabiliser
  of the edge $e_{\lambda_1}$.  Hence, for some word $f_2'$ in $H$
  \begin{eqnarray*}
    f_{k+1}\ r_{\lambda'_k} f_k \cdots r_{\lambda'_2} f_2\ 
    r_{\lambda'_1} f_1 h_1^{-1} h_1 & = &
    f_{k+1}\ r_{\lambda'_k} f_k \cdots r_{\lambda'_2} f'_2\ 
    r_{\lambda_1} h_1
  \end{eqnarray*}
  modulo $R_1$. By induction the two shorter h-products
  $h_{k+1}\ r_{\lambda_k} h_k \cdots r_{\lambda_2} h_2$ and 
  $f_{k+1}\ r_{\lambda'_k} f_k \cdots r_{\lambda'_2} f'_2$ are 
  equivalent modulo $R_0 \union R_1$, and so $p_1 = p_2$ modulo
  $R_0 \union R_1$.
\end{proof}

\begin{claim} \label{claim0+1+2}
  Suppose that two h-products represent the same element of $G$ and induce
  edge paths that are equivalent modulo backtracking.  Then they are
  equivalent modulo $R_0 \union R_1 \union R_2$.
\end{claim}
\begin{proof}
  It is enough to show that any h-product is equivalent to an h-product
  that represents a path without any backtracking.  Furthermore, if we
  proceed by induction on the length of the h-product, it is enough to
  show that any h-product whose associated path has backtracking at the
  end is equivalent to a shorter h-product.

  Suppose that
  $g = h_{k+3}\ r_{\lambda_{k+2}} h_{k+2}\ r_{\lambda_{k+1}} h_{k+1}\ g_k$
  is such an h-product, ie
  \begin{eqnarray*}
            v_k     & = & v_0 \cdot g_k \\
            v_{k+1} & = & v_0 \cdot r_{\lambda_{k+1}} h_{k+1}\ g_k \\
    v_{k+2}\ =\ v_k & = & v_0 \cdot r_{\lambda_{k+2}} h_{k+2}\ 
                            r_{\lambda_{k+1}} h_{k+1}\ g_k 
  \end{eqnarray*}
  and $g_k$ is a shorter h-product.  So, multiplying by
  $g_k^{-1} h_{k+1}^{-1}$, we find that 
  $r_{\lambda_{k+2}} h_{k+2}\ r_{\lambda_{k+1}}$ is an h-product with
  associated path $(v_0, v_0\cdot r_{\lambda_{k+1}}, v_0)$.  Suppose
  that $r_{\lambda'} h r_\lambda = h'$ is the $R_2$ relation corresponding
  to this path.  Then $\lambda = \lambda_{k+1}$ and
  $v_0 \cdot r_{\lambda'} h = v_0 \cdot r_{\lambda_{k+2}} h_{k+2}$. So
  $\lambda' = \lambda_{k+2}$ and $h_{k+2} h^{-1}$ is in the stabiliser
  of the edge $e_{\lambda_{k+1}}$.  Therefore there exists a word $f$
  in $H$ such that
  \[ 
    h_{k+3}\ r_{\lambda_{k+2}} h_{k+2}\ r_{\lambda_{k+1}} h_{k+1}\ g_k
    =  h_{k+3} f\ r_{\lambda'} h\ r_\lambda h_{k+1}\ g_k 
  \]
  modulo $R_1$.  Hence modulo $R_2$ this is equal to
  $h_{k+3} f h' h_{k+1} g_k$, a shorter h-product.
\end{proof}

\begin{claim} \label{claim0+1+2+3}
  Any h-product equal to the identity in $G$ is equivalent to the
  identity modulo $R_0 \union R_1 \union R_2 \union R_3$.
\end{claim}
\begin{proof}
  Given any h-product $g_k$ equal to the identity in $G$ its associated
  edge path must be a loop.  Since $X$ is simply-connected this loop is
  the boundary of a union of faces of $X$.  So choose one of these
  faces $f$ touching the loop at a vertex $v$ then modulo
  $R_0 \union R_1 \union R_2$ we can add backtracking starting at $v$
  going around the boundary of $f$.  Modulo $R_3$ we can remove one
  pass round $\boundary f$.  This leaves a new loop that can be spanned
  by one less face, which, by induction on the minimum number of
  faces needed to span a loop, is equivalent to the identity.
\end{proof}

\begin{proof}[Proof of Theorem \ref{mainthm}]
  Given any word in the generators, $S_0 \union S_1$, that is equal to
  the identity in $G$ then modulo $R_2$ it is equivalent to an h-product
  and so by Claim~\ref{claim0+1+2+3} is equivalent to the identity 
  modulo $R_0 \union R_1 \union R_2 \union R_3$.
\end{proof}

\section{The complex $\Xfull{n}$} \label{Xfull}

An embedded disc $D \subseteq \halfspace$ is said to \emph{cut out $a_i$}
if the interior of $D$ is disjoint from $a_*$, the arc $a_i$ is contained
in the boundary of $D$ and the boundary of $D$ lies in 
$a_i \union \boundary \halfspace$, ie $a_i \subset \boundary D$ and 
$\boundary D \subset a_i \union \boundary\halfspace$.
A \emph{cut system for $a_*$} is the isotopy class of $n$ pairwise
disjoint discs $\langle D_1,D_2,\ldots D_n\rangle$ where each $D_i$
cuts out the arc $a_i$. Say that two cut systems
$\langle D_1, D_2, \ldots, D_n \rangle$ and
$\langle E_1, E_2, \ldots, E_n \rangle$ differ by a simple $i$-move if
$D_i \intersect E_i = a_i$ and $D_j = E_j$ for all $j \neq i$.  If
this is the case we will suppress the non-changing discs and write
$\langle D_i \rangle \to \langle E_i \rangle$.

\begin{definition}
  Define the cut system complex $\Xfull{n}$ as follows.  The set of all
  cut systems for $a_*$ forms the vertex set $\Xfull{n}^0$.  Two vertices
  are connected by a single edge iff they differ by a simple move.
  Finally, glue faces into every loop of the following form, giving
  triangular and rectangular faces.
  \[
    \xymatrix @C\halfroottwo pt{
      \langle D_i\rangle\ar@{-}[rr]& &\langle D'_i\rangle\ar@{-}[dl]\\
                      & \langle D''_i \rangle\ar@{-}[ul] & }
    \qquad 
    \xymatrix{
      \langle D_i,D_j\rangle\ar@{-}[r]
                               & \langle D'_i,D_j\rangle\ar@{-}[d]  \\
      \langle D_i,D'_j\rangle\ar@{-}[u]
                               & \langle D'_i , D'_j \rangle\ar@{-}[l]
    }
  \]

  Define the basepoint to be $v_0 = \langle d_1, d_2, \ldots, d_n \rangle$ 
  where the $d_i$ are vertical discs below the $a_i$, see Figure~\ref{fig1}.
  \begin{figure}[!htb]
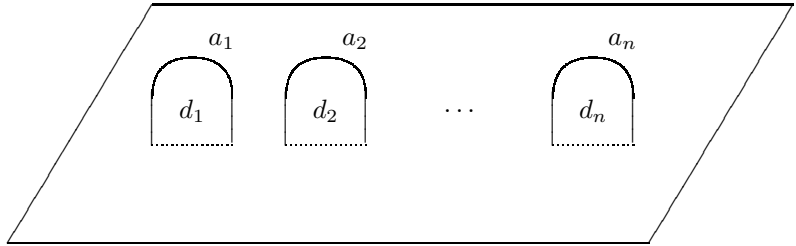

  \[ \xy
      <0pt,-50pt>;<20pt,-50pt>:
      (0.3,0.5);(12.3,0.5) **\dir{-}; 
      (3,5);(15,5) **\dir{-}; 
      (0.3,0.5);(3,5) **\dir{-}; 
      (12.3,0.5);(15,5) **\dir{-}; 
      (3,2.5); @(, (3,4)@+, (4.5,4)@+, (4.5,2.5), *\qspline{}, @i @); 
      (5.5,2.5); @(, (5.5,4)@+, (7,4)@+, (7,2.5)*\qspline{}, @i @); 
      (10.5,2.5); @(, (10.5,4)@+, (12,4)@+, (12,2.5)*\qspline{}, @i @); 
      (3,2.35); (4.5,2.35) **\dir{.}; 
      (5.5,2.35); (7,2.35) **\dir{.}; 
      (10.5,2.35); (12,2.35) **\dir{.}; 
      (8.75,3)*{\ldots};
      (4.3,4.3)*{a_1};
      (6.8,4.3)*{a_2};
      (11.8,4.3)*{a_n};
      (3.75,3)*{d_1};
      (6.25,3)*{d_2};
      (11.25,3)*{d_n};
  \endxy \]
  \caption{The arcs $a_i$ and the discs $d_i$} \label{fig1}
  \end{figure}
  Sometimes it is convenient to think of the $a_i$ and $d_i$ rotated
  by a quarter turn.
\end{definition}

Before we prove that this complex is simply connected we need the 
following lemma about substituting one disc for another.

Suppose that $v = \langle D_1, D_2, \ldots D_n \rangle$ is a vertex of
$\Xfull{n}$ with a choice of discs representing it and that $D$ and
$D^*$ are two discs cutting out the arc $a_i$.  We will say that the
tuple $(v, D, D^*)$ forms a \emph{valid substitution} if either $D\ne D_i$ or
$D=D_i$ and for all $j \ne i$, $D_j \intersect D^* = \emptyset$, ie if
$D$ is in $v$ then there exists an edge 
$\langle D = D_i \rangle \text{---} \langle D^* \rangle$.  If
$(v, D, D^*)$ forms a valid substitution then we can replace $D$ with
$D^*$ to get a vertex $v^*$, ie
\[ v^* = \begin{cases} v & \text{if $D_i \ne D$,} \\
               \langle D^*\rangle &\text{if $D_i = D$.} \end{cases} \]

Similarly, for any edge path $P$ with a choice of discs representing
each vertex, we say $(P, D, D^*)$ forms a valid substitution if for each
vertex $v$ of $P$ the tuple $(v, D, D^*)$ forms a valid substitution.  
If $(P, D, D^*)$ forms a valid substitution then we can replace each
occurrence of $D$ with $D^*$, ie replace each vertex $v$ with $v^*$,
giving a new path $P^*$.

\begin{lemma} \label{sub}
  If $(P, D, D^*)$ forms a valid substitution, where 
  $P = (v_1, \ldots, v_k)$, then $P^*$ is a path and the loop
  \[ \xymatrix{
              v_1   \ar@{-}[rr]^{P} \ar@{-}[d] & & v_k \ar@{-}[d] \\
              v_1^* \ar@{-}[rr]^{P^*}          & & v_k^*
  } \]
  is homotopic to a point. Moreover, if $P$ is a loop then so is $P^*$
  and they are homotopic as loops.
\end{lemma}
\begin{proof}
  Clearly we may assume that $D$ and $D^*$ are not isotopic, otherwise
  $P=P^*$.  Suppose that $D$ and $D^*$ cut out the arc $a_i$.  For each
  vertex $v$ of $P$ we have that either $v=v^*$ or $(v, v^*)$ is an edge of
  $\Xfull{n}$.

  For each edge $(u,v)$ in $P$, where $u=\langle D_j \rangle$ and
  $v = \langle D'_j \rangle$, we have the following possibilities.  If
  $D$ is not in $u$ nor $v$ then $(u,v) = (u^*, v^*)$.  Otherwise we
  have two cases depending on whether $i=j$ or not.

  If $i=j$ then only one of either $u$ or $v$ contains $D$.  Suppose that 
  $D \in u$, ie $D_j=D$.  If $D^*=D'_j$ then $u^*=v^*=v$ and $(u,v)$ is
  homotopic to $(u^*, v^*)$ in $\Xfull{n}^1$.  Otherwise, if $D^* \neq D_j$,
  we have the following face of $\Xfull{n}$.
  \[ \xymatrix{
              u \ar@{-}[r]^{P} \ar@{-}[d] & **[r] v = v^* \\
              u^*\ar@{-}[ur]_{P^*}
  } \]
  If $i\neq j$ the we have the following face of $\Xfull{n}$.
  \[ \xymatrix{
              \langle D, D_j \rangle \ar@{-}[r]^{P} \ar@{-}[d] 
                               & \langle D, D_j' \rangle \ar@{-}[d] \\
              \langle D^*, D_j \rangle \ar@{-}[r]^{P^*} 
                               & \langle D^*, D_j' \rangle
  } \]

  In either case there is a homotopy from $(u,v)$ to $(u^*,v^*)$ that
  agrees with the homotopies between the vertices of $P$ and $P^*$.
  Therefore $P$ is homotopic to $P^*$.
\end{proof}

\begin{theorem} \label{Xfull-simply-connected}
  The complex $\Xfull{n}$ is connected and simply connected.
\end{theorem}

\begin{proof}
  Given a loop in $\Xfull{n}$ it is homotopic to an edge path $P$.  Now
  choose discs to represent each vertex of $P$.  We shall write $D \in P$
  if $D$ is one of the discs chosen as a representative of some vertex
  of $P$. Assuming that the intersection of the discs $D \in P$ with
  $d_1 \union d_2 \union \ldots \union d_n$ isn't only
  $a_1, a_2, \ldots, a_n$.  We can carry out the following procedure.

  For some $i$ the union of the discs in $P$ intersects $d_i$ in a
  non-empty collection of arcs.  Pick an arc $\alpha$ of this intersection that
  is lowest in the sense that it doesn't separate the entirety of any
  other arc from $\boundary \halfspace \intersect d_i$. For example, see 
  Figure~\ref{arcs} where
  $\alpha$ and $\gamma$ are lowest but $\beta$ is not.
  \begin{figure}[!htb]
    \centering
    \begin{picture}(0,0)
      \put(10,34){$\alpha$}
      \put(70,29){$\beta$}
      \put(61,10){$\gamma$}
    \end{picture}
    \includegraphics{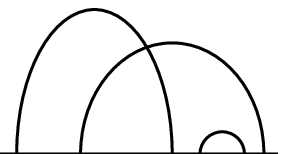}
    \caption{Lowest arcs $\alpha$ and $\gamma$} \label{arcs}
  \end{figure}

  The arc $\alpha$ comes from some $D \in P$.  Now cut $D$ along 
  $\alpha$, discard the section not containing $a_i$ and glue in a disc
  parallel to $d_i$.  This results in a new disc $D^*$ whose intersection
  with $d_i$ contains at least one less arc.  

  Any disc $E \in P$ for which $E \intersect D = a_j \text{ or } \emptyset$
  also has $E \intersect D^* = a_j \text{ or } \emptyset$ respectively;
  if not $E$ must intersect $D^*$ in the section parallel to $d_i$ and
  this contradicts the condition that $\alpha$ is a lowest arc.
  Therefore the triple $(P, D, D^*)$ form a valid substitution and, by
  Lemma~\ref{sub}, we can replace $D$ with $D^*$ to get a new homotopic
  loop $P^*$.

  We now have a homotopic loop $P^*$ that has fewer intersections with
  $d_1 \union d_2 \union \ldots \union d_n$.  So by induction on the
  number of intersections we have proved the following.

  \begin{claim*} The path $P$ is homotopic to a path whose vertices
    admit representative discs which intersect the discs
    $d_1, d_2, \ldots, d_n$ only in the arcs $a_1, a_2, \ldots, a_n$.
  \end{claim*}

  So we may assume that the path $P$ meets $d_1, d_2, \ldots, d_n$ only
  in the arcs $a_1, a_2, \ldots, a_n$.  Therefore, for each $D \in P$
  cutting out the arc $a_i$, $(P,D,d_i)$ forms a valid substitution and
  so by in turn replacing each $D\in P$ with $d_i$ we see that $P$ is
  homotopic to the constant path $v_0$.  The connectedness of $\Xfull{n}$
  follows by taking $P$ to be a constant loop.
\end{proof}

Up to homotopy the group $\Hilden{2n}$ acts on $(\halfspace,a_*)$ by 
homeomorphisms, therefore it takes cut systems to cut systems.  The
edges and faces of $\Xfull{n}$ are determined by the intersections of
pairs of discs, hence this action on $\Xfull{n}^0$ extends to a cellular
action on $\Xfull{n}$.

\begin{theorem}
  The action of $\Hilden{2n}$ on $\Xfull{n}^0$ is transitive.
\end{theorem}
\begin{proof}
  Given a vertex $\langle D_1, D_2, \ldots, D_n \rangle$ of $\Xfull{n}$,
  if we take each $i$ in turn and look at the intersection of $D_i$ with
  $\boundary\halfspace$.  We see that this defines a path from one end
  of $a_i$ to the other.  If we now move one end around this path until
  it is close to the other and then move it straight back to its starting
  point we have an element of $\Hilden{2n}$ that moves $D_i$ to $d_i$.
  Combining all of these we see that 
  $\langle D_1, D_2, \ldots, D_n \rangle$ is in the orbit of $v_0$, ie
  the action is transitive on $\Xfull{n}^0$.
\end{proof}

\section{The complex $\X{n}$} \label{X}

We now construct a subcomplex $\X{n}$ of $\Xfull{n}$ with the same 
vertex set but with fewer edges and faces. 

Given an edge $e = (\langle D \rangle ,\langle D' \rangle)$ of
$\Xfull{n}$ define its length, $l(e)$, to be the number of arcs
underneath $D \union D'$.  In other words, since
$\halfspace \setminus D \union D'$ has two components, one bounded and
one unbounded, we can define the length as follows
\[ 
  l(e) = \# \{ i \st a_i \text{ is contained in the bounded component of }
                     \halfspace \setminus D \union D' \}. 
\]
Given two edges $e$ and $e'$ with the same length there exists an element
of $\Hilden{2n}$ taking $e$ to $e'$.

Say a rectangle $(\langle D,E\rangle, \langle D',E\rangle, 
\langle D',E'\rangle, \langle D,E'\rangle)$ is \emph{nested} if 
$E\union E'$ lies in the bounded component of 
$\halfspace \setminus D \union D'$ or vice versa, ie if one pair of
changing discs lies underneath the other.

For $i \leq j$ let $\Torbit_{ij}$ denote the subcomplex consisting of
all triangular faces of $\Xfull{n}$ with shortest two edges of length
$i$ and $j$.  Note, this implies that the remaining edge has length $i+j$.
Given a rectangular face of $\Xfull{n}$ we have two cases depending on
whether it is nested or not. Let $\Rorbit_{ij}$ denote the subcomplex
consisting of all rectangular nested faces with inner edge of length $i$
and outer edge of length $j$.  For $i \leq j$ let $\Sorbit_{ij}$ denote
the subcomplex consisting of all non-nested rectangular faces with edges
of length $i$ and $j$.

\begin{definition} 
  Let $\X{n}$ be the subcomplex of $\Xfull{n}$ with the same vertex set,
  all edges of length 1 and 2 and all faces from $\Rorbit_{12}$,
  $\Sorbit_{11}$ and $\Torbit_{11}$, ie 
  $\X{n} = \Rorbit_{11} \union \Sorbit_{11} \union \Torbit_{11}$. As
  the length of an edge is invariant under the action of $\Hilden{2n}$
  on $\Xfull{n}$ this action preserves each $\Torbit_{ij}$, $\Rorbit_{ij}$
  and $\Sorbit_{ij}$ and so preserves $\X{n}$.
\end{definition}

\begin{figure}[!htb]
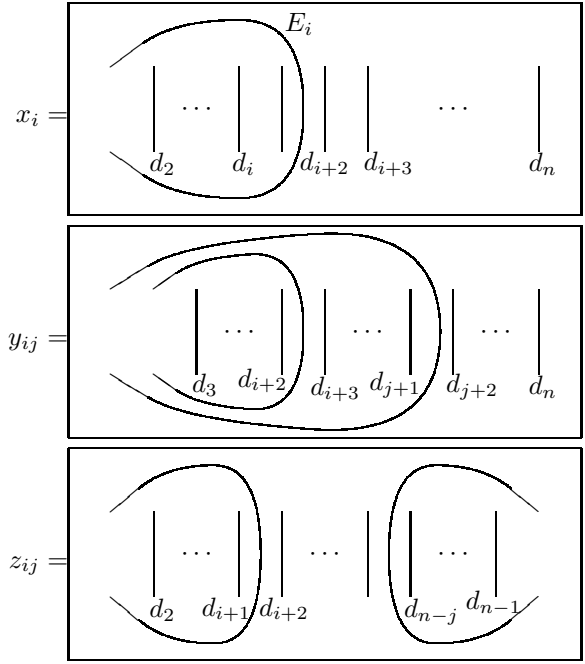

  \begin{align*} 
    x_i = & 
    \xy
      <0pt,-35pt>;<16pt,-35pt>:
      (0,0)*{}; (12,0)*{} **\dir{-}; 
      (0,5)*{}; (12,5)*{} **\dir{-}; 
      (0,0)*{}; (0,5)*{} **\dir{-}; 
      (12,0)*{}; (12,5)*{} **\dir{-}; 
      (2,1.5)*{}; (2,3.5)*{} **\dir{-}; 
      (3,2.5)*{\cdots};
      (4,1.5)*{}; (4,3.5)*{} **\dir{-}; 
      (5,1.5)*{}; (5,3.5)*{} **\dir{-}; 
      (6,1.5)*{}; (6,3.5)*{} **\dir{-}; 
      (7,1.5)*{}; (7,3.5)*{} **\dir{-}; 
      (9,2.5)*{\cdots};
      (11,1.5)*{}; (11,3.5)*{} **\dir{-}; 
      (1,1.5)*{}; @(, (2.5,0.4)@+, (5.5,0.4)@+, 
          (5.5,4.6)@+, (2.5,4.6)@+, (1,3.5)*\qspline{}, @i @); 
      (5.4,4.5)*{E_i};
      (2.2,1.2)*{d_2};
      (4.1,1.2)*{d_i};
      (6,1.2)*{d_{i+2}};
      (7.5,1.2)*{d_{i+3}};
      (11.1,1.2)*{d_n};
    \endxy \\
    y_{ij} = &
    \xy
      <0pt,-35pt>;<16pt,-35pt>:
      (0,0)*{}; (12,0)*{} **\dir{-}; 
      (0,5)*{}; (12,5)*{} **\dir{-}; 
      (0,0)*{}; (0,5)*{} **\dir{-}; 
      (12,0)*{}; (12,5)*{} **\dir{-}; 
      (3,1.5)*{}; (3,3.5)*{} **\dir{-}; 
      (4,2.5)*{\cdots};
      (5,1.5)*{}; (5,3.5)*{} **\dir{-}; 
      (6,1.5)*{}; (6,3.5)*{} **\dir{-}; 
      (7,2.5)*{\cdots};
      (8,1.5)*{}; (8,3.5)*{} **\dir{-}; 
      (9,1.5)*{}; (9,3.5)*{} **\dir{-}; 
      (10,2.5)*{\cdots};
      (11,1.5)*{}; (11,3.5)*{} **\dir{-}; 
      (1,1.5)*{}; @(, (2.75,0.4)@+, (8.7,0)@+, 
          (8.7,5)@+, (2.75,4.6)@+, (1,3.5)*\qspline{}, @i @); 
      (2,1.5)*{}; @(, (3.05,0.75)@+, (5.5,0.6)@+,
          (5.5,4.4)@+, (3.05,4.25)@+, (2,3.5)*\qspline{}, @i @); 
      %
      %
      (3.2,1.25)*{d_3};
      (4.55,1.35)*{d_{i+2}};
      (6.25,1.15)*{d_{i+3}};
      (7.65,1.25)*{d_{j+1}};
      (9.45,1.2)*{d_{j+2}};
      (11.1,1.2)*{d_n};
    \endxy \\
 z_{ij} = &
    \xy
      <0pt,-35pt>;<16pt,-35pt>:
      (0,0)*{}; (12,0)*{} **\dir{-}; 
      (0,5)*{}; (12,5)*{} **\dir{-}; 
      (0,0)*{}; (0,5)*{} **\dir{-}; 
      (12,0)*{}; (12,5)*{} **\dir{-}; 
      (2,1.5)*{}; (2,3.5)*{} **\dir{-}; 
      (3,2.5)*{\cdots};
      (4,1.5)*{}; (4,3.5)*{} **\dir{-}; 
      (5,1.5)*{}; (5,3.5)*{} **\dir{-}; 
      (6,2.5)*{\cdots};
      (7,1.5)*{}; (7,3.5)*{} **\dir{-}; 
      (8,1.5)*{}; (8,3.5)*{} **\dir{-}; 
      (9,2.5)*{\cdots};
      (10,1.5)*{}; (10,3.5)*{} **\dir{-}; 
      (1,1.5)*{}; @(, (2.3,0.4)@+, (4.5,0.4)@+, 
          (4.5,4.6)@+, (2.3,4.6)@+, (1,3.5)*\qspline{}, @i @); 
      (11,1.5)*{}; @(, (9.7,0.4)@+, (7.5,0.4)@+, 
          (7.5,4.6)@+, (9.7,4.6)@+, (11,3.5)*\qspline{}, @i @); 
      %
      %
      (2.2,1.2)*{d_2};
      (3.7,1.2)*{d_{i+1}};
      (5.05,1.2)*{d_{i+2}};
      (8.5,1.15)*{d_{n-j}};
      (9.95,1.35)*{d_{n-1}};
    \endxy 
  \end{align*}
  \caption{The vertices $x_i$, $y_{ij}$ and $z_{ij}$}\label{fig2}
\end{figure}

A vertex $v = \langle D_1, \ldots, D_n \rangle$ is completely determined
by the intersection of the discs $D_i$ with $\boundary \halfspace$.
Using this we can define the vertices $x_i$ for $0\le i\le n-1$, $y_{ij}$
for $0\le i\le n-2$ and $j=0$ or $i<j\le n-1$ and $z_{ij}$ for 
$0\le i, j, i+j \le n-2$ as in Figure~\ref{fig2}. So we have 
$l(v_0, x_i) = i$, $l(v_0, y_{0j}) = j$, $l(v_0, y_{i0}) = i$,
$l(v_0, z_{i0}) = i$ and $l(v_0, z_{0j}) = j$.  Note, there is some 
redundancy in this notation, ie $x_i = y_{0i}$ and 
$x_0 = y_{00} = z_{00} = v_0$.

We now define the faces $R_{ij} \in \Rorbit_{ij}$, $S_{ij} \in \Sorbit_{ij}$,
$T_{ij} \in \Torbit_{ij}$ of $\Xfull{n}$ as follows.
\[
  R_{ij}\ = \vcenter{
              \xymatrix{ 
                   y_{00} \ar @{-} [r] \ar@{-}[d] & y_{i0}\ar@{-}[d]\\
                   y_{0j} \ar@{-}[r]              & y_{ij} } } \qquad
  S_{ij}\ = \vcenter{ 
              \xymatrix{ 
                   z_{00} \ar @{-} [r] \ar@{-}[d] & z_{i0}\ar@{-}[d]\\
                   z_{0j} \ar@{-}[r]              & z_{ij} } } \qquad
  T_{ij}\ =\! \vcenter{
              \xymatrix @C\halfroottwo pt{ 
                   x_0 \ar@{-}[rr] & & x_i \ar@{-}[dl] \\
                         & x_{i+j} \ar@{-}[ul] & } }
\]
For every face in $\Xfull{n}$ there is an element of $\Hilden{2n}$
taking it to one of these representatives.

\begin{theorem} \label{X-simply-connected}
  The complex $\X{n}$ is simply connected.
\end{theorem}

\begin{proof}
  Figure~\ref{fig3} shows that the boundary of each of the faces $R_{ij}$
  for $1 < i < j$ and $S_{ij}$ for $1 < i,j$ can be expressed as the boundary
  of a union of faces with shorter edges. The first column shows how to
  replace faces where the first index is not 1.  Then the second column
  can be used to reduce the second index to either 2 or 1 respectively.

  \begin{figure}[!htb]
    \[ \begin{array}{cc}
      R_{ij} \colon \!\!\! \vcenter{ 
        \xymatrix{
           y_{00} \ar@{-}[rr]|-{}="0" \ar@{-}[ddd]|-{}="1" \ar@{-}[dr] 
                               & & y_{i0} \ar@{-}[ddd]|-{}="2" \ar@{-}[dl] \\
           & y_{i-1,0} \ar@{-}[d]|-{}="3" & \\
           & y_{i-1,j} \ar@{-}[dl] \ar@{-}[dr] & \\
           y_{0j} \ar@{-}[rr]|-{}="4" & & y_{ij} 
            \ar@{}"0";"2,2"|-{\Torbit_{1,i-1}}
            \ar@{}"4";"3,2"|-{\Torbit_{1,i-1}}
            \ar@{}"1";"3"|-{\Rorbit_{i-1,j}}
            \ar@{}"2";"3"|-{\Rorbit_{1j}}
        }
      }
     &
      R_{1j} \colon \!\!\! \vcenter{ 
        \xymatrix{
           y_{00} \ar@{-}[rrr]|-{}="0" \ar@{-}[dr] \ar@{-}[dd]|-{}="1" 
             & & & y_{10} \ar@{-}[dl] \ar@{-}[dd]|-{}="3" \\
          & y_{0,j-1} \ar@{-}[r]|-{}="2" \ar@{-}[dl]
             & y_{1,j-1} \ar@{-}[dr] & \\
           y_{0j} \ar@{-}[rrr]|-{}="4" & & & y_{1j}
             \ar@{}"0";"2"|-{\Rorbit_{1,j-1}}
             \ar@{}"1";"2,2"|-{\Torbit_{1,j-1}}
             \ar@{}"2,3";"3"|-{\Torbit_{1,j-1}}
             \ar@{}"2";"4"|-{\Sorbit_{11}}
        }
       }
     \\
      S_{ij} \colon \!\!\! \vcenter{ 
       \xymatrix{
           z_{00} \ar@{-}[rr]|-{}="1" \ar@{-}[ddd]|-{}="2" \ar@{-}[dr] 
                                 & & z_{i0} \ar@{-}[ddd]|-{}="4" \ar@{-}[dl] \\
           & z_{i-1,0} \ar@{-}[d]|-{}="3" & \\
           & z_{i-1,j} \ar@{-}[dl] \ar@{-}[dr] & \\
           z_{0j} \ar@{-}[rr]|-{}="5" & & z_{ij} 
              \ar@{}"1";"2,2"|-{\Torbit_{1,i-2}}
              \ar@{}"2";"3"|-{\Sorbit_{i-1,j}}
              \ar@{}"3";"4"|-{\Sorbit_{1j}}
              \ar@{}"5";"3,2"|-{\Torbit_{1,i-2}}
        }
      }
     &
      S_{1j} \colon \!\!\! \vcenter{ 
        \xymatrix{ \labelstyle{\displaystyle} \objectstyle{\displaystyle}
           z_{00} \ar@{-}[rrr]|-{}="1" \ar@{-}[dr] \ar@{-}[dd]|-{}="2" 
              & & & z_{10} \ar@{-}[dl] \ar@{-}[dd]|-{}="4" \\
          & z_{0,j-1} \ar@{-}[r]|-{}="3" \ar@{-}[dl] 
              & z_{1,j-1} \ar@{-}[dr] &  \\
           z_{0j} \ar@{-}[rrr]|-{}="5"  & & & z_{1j}
              \ar@{}"1";"3"|-{\Sorbit_{1,j-1}}
              \ar@{}"2";"2,2"|-{\Torbit_{1,j-1}}
              \ar@{}"2,3";"4"|-{\Torbit_{1,j-1}}
              \ar@{}"3";"5"|-{\Sorbit_{11}}
        }
       }
    \end{array} \]
    \caption{Decomposing rectangular faces}\label{fig3}
  \end{figure}
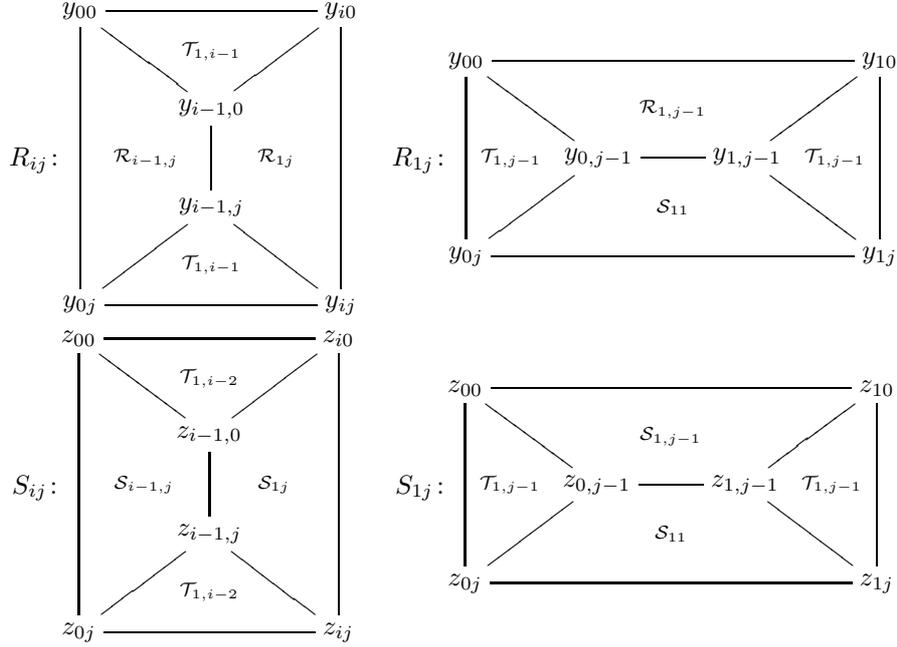

  As each of the rectangular faces can be moved to one of $R_{ij}$ or
  $S_{ij}$ by some element of $\Hilden{2n}$ it follows that every loop
  in $\X{n}$ is null-homotopic in 
  \[ \Rorbit_{12} \union \Sorbit_{11} 
                \union \Union_{1 \leq i \leq j \leq n}\Torbit_{ij}. \]

  Let the $E_i$ be the discs as shown in Figure~\ref{fig2}, ie
  $x_i = \langle E_i, d_2, d_3, \ldots, d_n \rangle$.  For $j > 2$ let
  $A_j$ the be full subcomplex of $\Xfull{n}$ containing all the vertices
  ``between'' $x_0$ and $x_j$, ie
  \begin{multline*} 
    A_j^0 = \{ \langle D, d_2, d_3, \ldots, d_n \rangle \in \Xfull{n}^0 \\
            \st D \ne d_0 \text{ or } E_j, 
                \text{interior of $D \subset$ bounded component of 
                    $\halfspace \setminus E_0 \union E_j$} \}.
  \end{multline*}
  Choose $x_1$ as a base point of $A_j$.

  For every edge $(u, v)$ of $A_j$ we have the following two triangles
  in $\Xfull{n}$.  Note that all edges have length less than $j$.
  \[ \xymatrix{ 
                    & \ar@{-}[dl] x_0 \ar@{-}[dr] &                 \\
        u \ar@{-}[rr] \ar@{-}[dr] & & v \ar@{-}[dl]                 \\
                                & x_j & 
  } \]

  \begin{lemma}\label{between}
    The subcomplex $A_j$ is path connected.
  \end{lemma}
  \begin{proof}
    Given a vertex $v = \langle D, d_2, \ldots d_n \rangle \in A_j^0$.
    First suppose that for some $2 < i \leq j$ there exists a path
    $\gamma$ on $\boundary \halfspace$ from $d_i$ to $E_j$ such that
    $\gamma$ does not cross $E_1$, $D$ or $d_l$ for $l\neq i$.  Let 
    $D'$ be a disc parallel to $E_j$ except in a neighbourhood of 
    $\gamma$ where we glue in the boundary of a neighbourhood of
    $\gamma \union d_i$.  Then there is a path $(v, v', x_1)$ in $A_j$
    where $v' = \langle D' \rangle$.  See  Figure~\ref{tunnel}.

    \begin{figure}[!htb]
      \centering
      \begin{picture}(0,0)
        \put(47,20){$E_1$}
        \put(84,75){$i$}
        \put(94,63){$\gamma$}
        \put(115,40){$D$}
        \put(200,0){$E_j$}
        \put(121,84){{\red $D'$}}
      \end{picture}
      \includegraphics{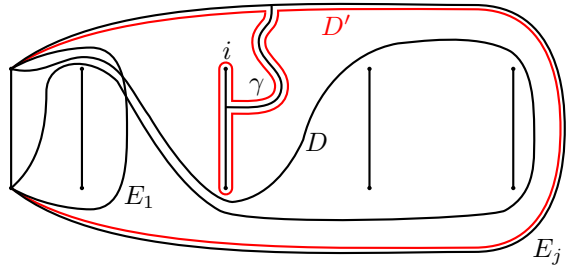}
      \caption{Tunnelling along $\gamma$} \label{tunnel}
    \end{figure}

    Now suppose that no such path exists on $\boundary \halfspace$.
    Each vertex $u = \langle D_u \rangle$ of $A_j$
    partitions the set $\{ d_2, d_3, \ldots d_{j+1} \}$ into two 
    non-empty subsets.  The first containing those discs that are
    between $d_1$ and $D_u$, the second those between $D_u$ and $E_j$.
    (If one of these sets were empty then we would have that either 
    $D_u = d_1$ or $D_u = E_j$.)  As $j>2$ at least one of these sets
    contains more than one disc.  Choose an $i\ne 1$ such that $d_i$
    is in this set.

    Now draw a path $\gamma$ on $\boundary \halfspace$ from $d_i$ to 
    $E_j$ that doesn't intersect $d_l$ for $l = 3, \ldots, j$ or $E_1$
    and only intersects $D$ transversely.  Starting at $d_i$ move along
    $\gamma$ and label the successive points of $\gamma \intersect D$
    as $p_1, p_2, \ldots, p_k$.  Now we can construct a sequence of discs
    $D = D^0, D^1, \ldots, D^k$ where each $D^{l+1}$ is parallel to $D^l$
    except in a neighbourhood of $p_{l+1}$ where we glue in the boundary
    of a sufficiently small neighbourhood of the disc $d_i$ and the 
    segment of $\gamma$ up to $p_{l+1}$.  With each successive $D^l$ the
    disc $d_i$ moves from one side of the partition to the other.  At 
    each step neither side of the partition is empty so 
    $\langle D^l \rangle$ is a vertex of $A_j$. This gives a path 
    $(v = \langle D^0\rangle, \langle D^1\rangle, \ldots, 
    \langle D^k \rangle)$ in $A_j$.  Now, $\langle D^k \rangle$
    satisfies the hypothesis above, therefore this path can be continued
    to the base point $x_1$.
  \end{proof}

  We can now complete the proof of Theorem~\ref{X-simply-connected}.
  So far we have shown that any loop in $\X{n}$ is the boundary of a
  union of faces in $\Rorbit_{12} \union \Sorbit_{11} \union
  \Union_{1\leq i\leq j\leq n}\Torbit_{ij}$.  For a given loop take an
  edge $(u,v)$ of maximal length $j$ in this union.  If $j > 2$ then the
  faces on either side of $(u,v)$ must be triangular with the remaining
  edges of length less than $j$.  So we have the following situation
  for some $u', v' \in \Xfull{n}^0$.
  \[ \xymatrix{ 
         & \ar@{-}[dl] u \ar@{-}[dd] \ar@{-}[dr] &      \\
       u' \ar@{-}[dr] & & v' \ar@{-}[dl]                \\
                     & v &                              } \]
  By Lemma~\ref{between} we can replace these two triangles with the
  following.
  \[ \xymatrix{ 
       &  & \ar@{-}[dll] \ar@{-}[dl] u \ar@{-}[dr] \ar@{-}[drr] & & \\
    u_0 \ar@{-}[drr] \ar@{-}[r] & 
    u_1 \ar@{-}[dr] \ar@{-}[r] & 
    \cdots \ar@{-}[r] & 
    u_{k-1} \ar@{-}[dl] \ar@{-}[r] & u_k \ar@{-}[dll] \\
                  & & v & &                             } \]
  Where $u_0 = u'$ and $u_k = v'$.  Here each edge has length less than
  $j$.  Therefore all edges of length greater that 2 can be replaced
  and so the loop is null-homotopic in $\X{n}$.
\end{proof}

\section{Calculating the presentation} \label{calc}

By Section~\ref{X} we have an $\Hilden{2n}$--action on a simply connected
cellular complex.  So we can now follow the method given in 
Section~\ref{method}.  

Using the fact that $\Hilden{2n}$ is a subgroup of $\Braid{2n}$, we can
define the following elements of $\Hilden{2n}$ in terms of
$\sigma_1, \ldots, \sigma_{2n-1}$ the generators of $\Braid{2n}$.
\begin{align*}
   r_1 & = \sigma_2 \sigma_1 \sigma_3^{-1} \sigma_2^{-1} && \\
   r_2 & = \sigma_4 \sigma_3 \sigma_2 \sigma_1 
           \sigma_5^{-1} \sigma_4^{-1} \sigma_3^{-1} \sigma_2^{-1}&& \\
   s_i & = \sigma_{2i} \sigma_{2i-1} \sigma_{2i+1} \sigma_{2i} 
               && \for i \in \{1, \ldots, n-1\}\\
   t_i & = \sigma_{2i-1}
               && \for i \in \{1, \ldots, n\}
\end{align*}
So $r_1$ is the first arc passing through the second, $r_2$ is the first
two arcs passing through the third, $s_i$ is the $i$th and $i+1$st arcs 
crossing and $t_i$ is the $i$th arc performing a half twist.  Subsequently
we will prove that these generate $\Hilden{2n}$.

\begin{proposition}
  The stabiliser of the vertex $v_0$ is isomorphic to the framed braid
  group and hence has a presentation $\langle S_0 \mid R_0 \rangle$ where
  \begin{gather*}
    \begin{array}{rclr@{\ }c@{\ }ll}
    S_0 &=& \multicolumn{5}{l}{\{ s_1, s_2, \ldots, s_{n-1},\ 
                                  t_1, t_2, \ldots, t_n\}}  \\
    R_0 &=& \big\{ &     s_i s_j &=& s_j s_i     & \for |i-j|>1, \\
        & &        & s_i s_j s_i &=& s_j s_i s_j & \for |i-j|=1, \\
        & &        &     t_i t_j &=& t_j t_i     & \text{for all } i, j,\\
        & &        &     s_i t_j &=& t_j s_i     
                                    & \text{if } j \notin \{i, i+1\},\\
        & &        &     s_i t_j &=& t_k s_i     
                                    & \text{if } \{i, i+1\} = \{j, k\} 
           \quad\big\}
        \end{array}
  \end{gather*}
\end{proposition}
\begin{proof}
  If we restrict to $\boundary \halfspace$, elements of $\Hilden{2n}$
  can be thought of as motions of the end points of the $a_i$.  For
  elements of the vertex stabiliser this motion moves the
  $d_i \intersect \boundary \halfspace$ among themselves, ie this is
  the fundamental group of configurations of $n$ line segments in the
  plain, the framed braid group.
\end{proof}

We have two edge orbits, one consisting of edges of length 1 and the
other consisting of edges of length 2.  Note that our choice of $r_1$
and $r_2$ mean that
\begin{gather*}
 (v_0, v_0 \cdot r_1) \in l^{-1}(1) \\
 (v_0, v_0 \cdot r_2) \in l^{-1}(2).
\end{gather*}
For $i=1,\ 2$, let $I_i$ denote the stabiliser of the edge
$(v_0, v_0 \cdot r_i)$, ie the subgroup of all elements that fix both
$v_0$ and $v_0 \cdot r_i$.

\begin{proposition}
  The subgroups $I_1$ and $I_2$ are generated as follows.
  \begin{align*}
    I_1 & = \langle t_2, t_3, \ldots, t_n,\ s_3, s_4,\ldots, s_{n-1},\ 
                 s_1 s_1 t_1 t_1,\ s_2 s_1 s_1 s_2 \rangle \\
    I_2 & = \langle t_2, t_3,\ldots, t_n,\ s_2,\ s_4, s_5,\ldots, s_{n-1},\  
        s_1 s_2 s_2 s_1 t_1 t_1,\ s_3 s_2 s_1 s_1 s_2 s_3 \rangle
  \end{align*}
\end{proposition}
\begin{proof}
  For $I_1$ [$I_2$] the motion of the $d_i$ outside of $d_1 \union E_2$
  [$d_1 \union E_3$] is generated by $t_3, t_4, \ldots, t_n, s_3, s_4,
  \ldots, s_{n-1}$ and $s_2 s_1 s_1 s_2$ [$t_4, t_5, \ldots, t_n$,
  $s_4, s_5, \ldots, s_{n-1}$ and $s_3 s_2 s_1 s_1 s_2 s_3$], the motion
  of the $d_i$ inside $d_1 \union E_2$ [$d_1 \union E_3$] is generated
  by $t_2$ [$t_2, t_3, s_2$] and the motion of $d_1 \union E_2$
  [$d_1 \union E_3$] is generated by $s_1 s_1 t_1 t_1$
  [$s_1 s_2 s_2 s_1 t_1 t_1$].
\end{proof}

We are now ready to calculate relations for $R_1$, $R_2$ and $R_3$.
The following relations are easily verifiable, in fact most of them
take place in $\Braid{8}$.

\subsection*{The $R_1$ relations}
To calculate the $R_1$ relations we have to find, for each edge orbit
representative $(v_0, v_0 \cdot r_i)$ and each generator $t$ of $I_i$,
a word $h$ in $S_0$ such that $r_i t r_i^{-1} = h$.  One possibility
is the following.
\begin{align}
    r_1 t_2 r_1^{-1} & = t_1                       \tag{$R_1$1}  \label{1-R1} \\
    r_1 t_k r_1^{-1} & = t_k & \text{for } k > 2   \tag{$R_1$2}  \label{2-R1}  \\
    r_1 s_k r_1^{-1} & = s_k & \text{for } k > 2   \tag{$R_1$3}  \label{3-R1}  \\
    r_1 s_1 s_1 t_1 t_1 r_1^{-1} & = s_1 s_1 t_2 t_2 &
                                                   \tag{$R_1$4}  \label{4-R1}  \\
    r_1 s_2 s_1 s_1 s_2 r_1^{-1} & = s_2 s_1 s_1 s_2 &
                                                   \tag{$R_1$5}  \label{5-R1}  \\[1em]
    r_2 t_2 r_2^{-1} & = t_1 &                     \tag{$R_1$6}  \label{6-R1}  \\
    r_2 t_3 r_2^{-1} & = t_2 &                     \tag{$R_1$7}  \label{7-R1}  \\
    r_2 t_k r_2^{-1} & = t_k & \text{for } k > 3   \tag{$R_1$8}  \label{8-R1}  \\
    r_2 s_2 r_2^{-1} & = s_1 &                     \tag{$R_1$9}  \label{9-R1}  \\
    r_2 s_k r_2^{-1} & = s_k & \text{for } k > 3   \tag{$R_1$10} \label{10-R1} \\
    r_2 s_1 s_2 s_2 s_1 t_1 t_1 r_2^{-1} & = s_2 s_1 s_1 s_2 t_3 t_3 &
                                                   \tag{$R_1$11} \label{11-R1} \\
    r_2 s_3 s_2 s_1 s_1 s_2 s_3 r_2^{-1} & = s_3 s_2 s_1 s_1 s_2 s_3 &
                                                   \tag{$R_1$12} \label{12-R1}
\end{align}

\subsection*{The $R_2$ relations}
To calculate the $R_2$ relations we need to find, for each edge orbit
representative $(v_0, v_0 \cdot r_i)$, an h-product $r_i h\ r_i$ for
the path $(v_0, v_0 \cdot r_i, v_0)$ and a word $h'$ in $S_0$ such that
$r_i h\ r_i = h'$.
\begin{align}
    r_1 t_1 s_1\ r_1 & = s_1 t_1                   \tag{$R_2$1}  \label{1-R2} \\
    r_2 s_1 t_2 s_2\ r_2 & = s_2 s_1 t_1           \tag{$R_2$2}  \label{2-R2}
\end{align}

\subsection*{The $R_3$ relations}
To calculate the $R_3$ relations we need to find, for each edge orbit,
an h-product representing the boundary of a face in the orbit and an
equivalent word in $S_0$.  The following are such relations for the
$\Sorbit_{11}$, $\Rorbit_{12}$ and $\Torbit_{11}$ orbits respectively.
\begin{gather}
  \begin{split}
     r_1 s_1 s_2 s_3 s_1 s_2\ r_1 s_1 s_2 s_3 s_1 s_2 t_2 t_4\ 
     r_1 s_2 s_3 s_1 s_2\ r_1 \\
              = s_1 s_2 s_3 s_1 s_2 s_1 s_2 s_1 s_3 s_2 s_2 s_3 s_1 s_2 t_1 t_3
  \end{split}                                      \tag{$R_3$1}  \label{1-R3} \\
     r_1\ r_2 s_1 s_2 s_1 t_2 t_3\ r_1\ r_2 
              = s_2 s_1 s_2 t_1 t_2                \tag{$R_3$2}  \label{2-R3} \\
     r_2 s_1 t_2\ r_1 s_2 s_1\ r_1 = s_1 s_2 s_1 t_1
                                                   \tag{$R_3$3}  \label{3-R3}
\end{gather}

\begin{figure}[!htb]
  \[ \xymatrix{ 
    \includegraphics[scale=0.09]{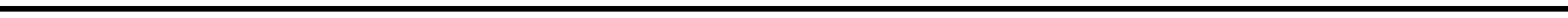} \ar@{->}[r]  
          & \includegraphics[scale=0.09]{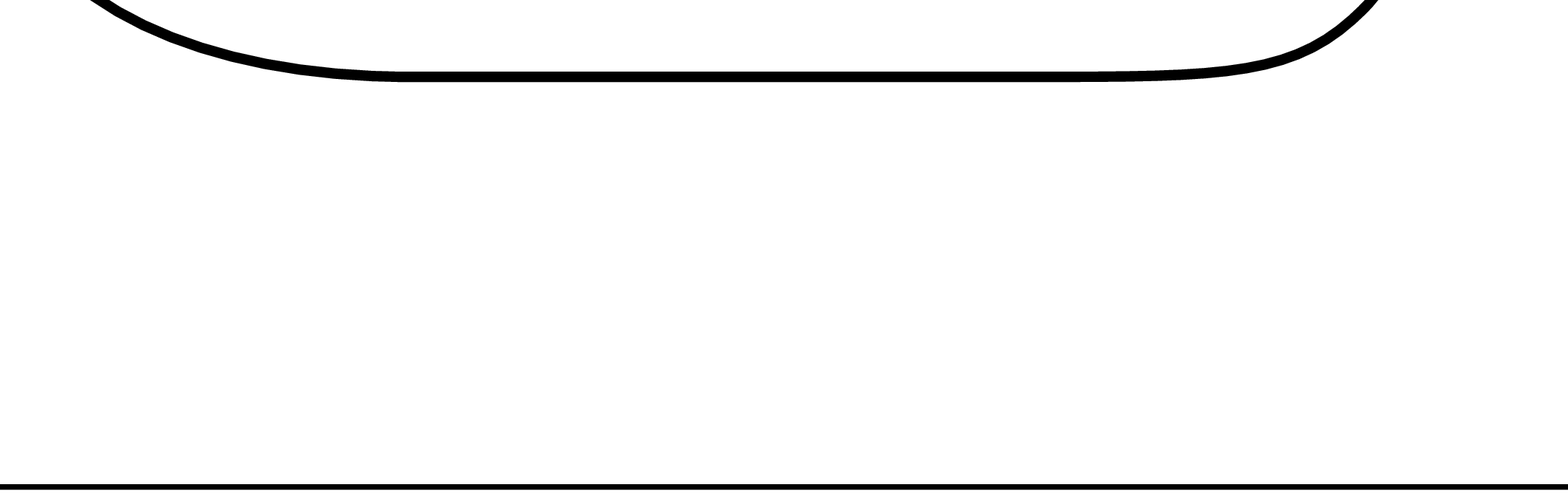} \ar@{->}[d] \\
    \includegraphics[scale=0.09]{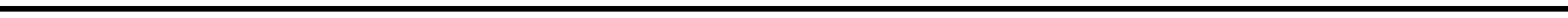} \ar@{->}[u] 
          & \includegraphics[scale=0.09]{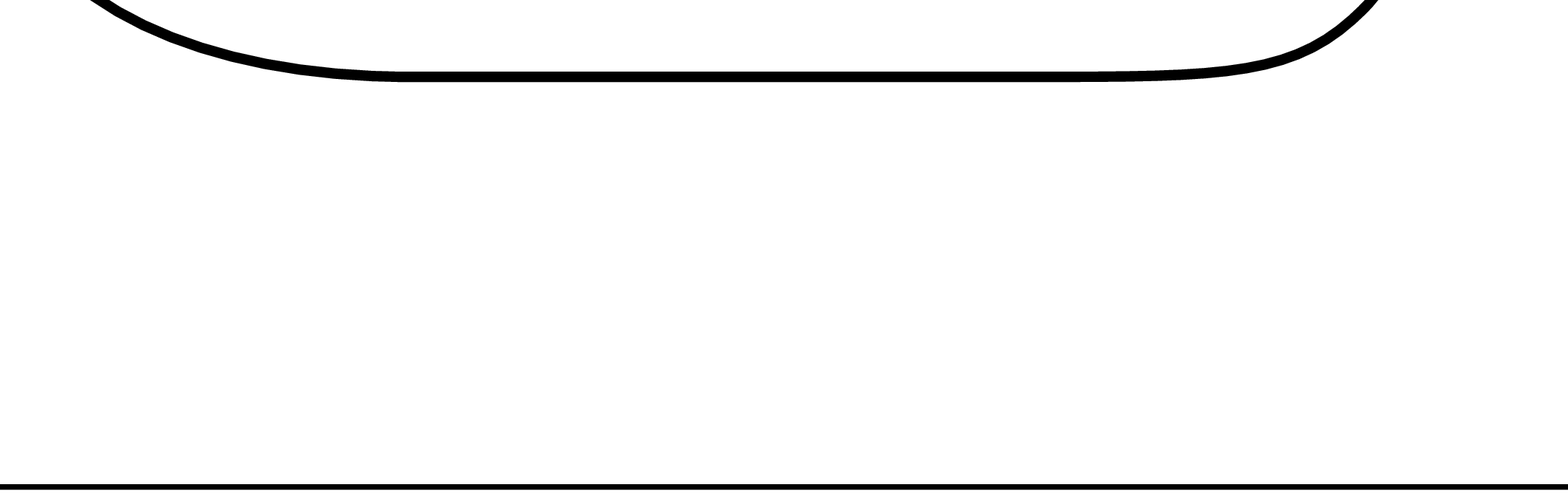} \ar@{->}[l]
   } \]
  \caption{The path given by the h-product on the LHS of \eqref{1-R3}}
  \label{1-R3:diagram}
\end{figure}

If we use a different set of generators, similar to those found by
Hilden, then we can get a more braid like presentation.  Let 
$p_i = \sigma_{2i} \sigma_{2i-1} \sigma_{2i+1}^{-1} \sigma_{2i}^{-1}$
for $1 \leq i < n$.  So $p_i$ is the $i$th arc passing under the $i+1$st
arc.

\begin{theorem}
  The group $\Hilden{2n}$ has a presentation with generators $p_i$,
  $s_j$ and $t_k$ for $1 \leq i,j<n$ and $1 \leq k \leq n$ and the 
  following relations.
  \begin{align}
              p_i p_j & = p_j p_i         & \text{for } |i-j|>1
                                                        \tag{P1}  \label{P1}  \\
          p_i p_j p_i & = p_j p_i p_j     & \text{for } |i-j|=1
                                                        \tag{P2}  \label{P2}  \\
              s_i s_j & = s_j s_i         & \text{for } |i-j|>1
                                                        \tag{P3}  \label{P3}  \\
          s_i s_j s_i & = s_j s_i s_j     & \text{for } |i-j|=1
                                                        \tag{P4}  \label{P4}  \\
              p_i s_j & = s_j p_i         & \text{for } |i-j|>1
                                                        \tag{P5}  \label{P5}  \\
      p_i s_{i+1} s_i & = s_{i+1} s_i p_{i+1}           \tag{P6}  \label{P6}  \\
  p_{i+1} p_i s_{i+1} & = s_i p_{i+1} p_i               \tag{P7}  \label{P7}  \\
  p_{i+1} s_i s_{i+1} & = s_i s_{i+1} p_i               \tag{P8}  \label{P8}  \\
      p_i t_i s_i p_i & = s_i t_i                       \tag{P9}  \label{P9}  \\
              p_i t_j & = t_j p_i   & \text{for } j \neq i, \text{ or } i+1
                                                        \tag{P10} \label{P10} \\
          p_i t_{i+1} & = t_i p_i                       \tag{P11} \label{P11} \\
              s_i t_j & = t_j s_i   & \text{if } j \neq i \text{ or } i+1
                                                        \tag{P12} \label{P12} \\
              s_i t_j & = t_k s_i   & \text{if } \{i, i+1\} = \{j, k\}
                                                        \tag{P13} \label{P13} \\
              t_i t_j & = t_j t_i   & \text{for } 1 \leq i,j \leq n
                                                        \tag{P14} \label{P14}
  \end{align}
\end{theorem}

\begin{figure}[!htb]
  \centering
  \begin{tabular}{ccc}
    \begin{minipage}[c]{3em}
      \centering \includegraphics{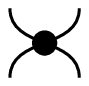}
    \end{minipage} & \begin{minipage}[c]{3em}
      \centering \includegraphics{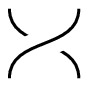} 
    \end{minipage} & \begin{minipage}[c]{3em}
      \centering \includegraphics{images/t.eps} 
    \end{minipage} \\
    $p$ & $s$ & $t$ \\
    \begin{minipage}[c]{3em}
      \centering \includegraphics{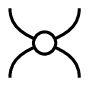}
    \end{minipage} & \begin{minipage}[c]{3em}
      \centering \includegraphics{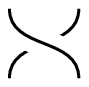} 
    \end{minipage} & \begin{minipage}[c]{3em}
      \centering \includegraphics{images/T.eps} 
    \end{minipage} \\
    $p^{-1}$ & $s^{-1}$ & $t^{-1}$ 
  \end{tabular}
  \caption{Pictorial representation of the $p$, $s$, $t$ and their inverses}
  \label{pst}
\end{figure}

These generators can be represented pictorially as in Figure~\ref{pst}.

\begin{figure}[!htb]
  \centering
  $\begin{minipage}[c]{6em}
      \hfill \includegraphics{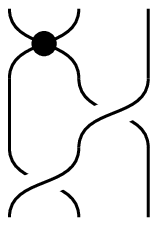}
    \end{minipage} = \begin{minipage}[c]{6em}
      \includegraphics{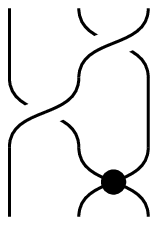}
    \end{minipage}$ 
    $\begin{minipage}[c]{4em}
      \hfill \includegraphics{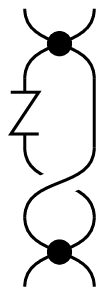}
    \end{minipage} = \begin{minipage}[c]{4em}
      \includegraphics{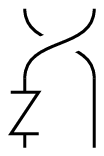}
    \end{minipage}$
  \caption{Pictorial representation of \eqref{P6} and \eqref{P9}}
  \label{relations}
\end{figure}

\begin{proof}
  Since $\Hilden{2n}$ is a subgroup of the braid group it is easy to
  check that these relations all hold.  So it remains to prove that
  each of the relations in $R_0$, $R_1$, $R_2$ and $R_3$ can be deduced
  from \eqref{P1}--\eqref{P14} using the fact that $r_1 = p_1$ and
  $r_2 = p_2 p_1$.  First note that $R_0$ is a subset of these relations.
  The relations \eqref{1-R1}, \eqref{2-R1}, \eqref{3-R1} and \eqref{1-R2}
  follow directly from \eqref{P11}, \eqref{P10}, \eqref{P5} and \eqref{P9}
  respectively.  The remaining relations can be deduced as follows.

  \newcommand{\deductionTableWidth}{xxxxxx & xxxxxxxxxxxxxx & xxxxxxxxxxxxxxxxxxxxxx & xxxxxxxxxx \kill}

  \begin{longtable}{cr@{$\ =\ $}lc} \deductionTableWidth
    \eqref{4-R1}: & $r_1 s_1 s_1 t_1 t_1 r_1^{-1}$ 
       & $p_1 \underline{s_1 s_1 t_1} t_1 p_1^{-1}$      & $\eqref{P13}^2$   \\
     & & $\underline{p_1 t_1 s_1} s_1 t_1 p_1^{-1}$      & $\eqref{P9}$      \\
     & & $s_1 t_1 \underline{p_1^{-1} s_1 t_1 p_1^{-1}}$ & $\eqref{P9}$      \\
     & & $s_1 \underline{t_1 t_1 s_1}$                   & $\eqref{P13}^2$   \\
     & & $s_1 s_1 t_2 t_2$                                                   \\
  \end{longtable}

  \begin{longtable}{cr@{$\ =\ $}lc} \deductionTableWidth
    \eqref{5-R1}: & $r_1 s_2 s_1 s_1 s_2 r_1^{-1}$
       & $\underline{p_1 s_2 s_1} s_1 s_2 p_1^{-1}$      & $\eqref{P6}$      \\
     & & $s_2 s_1 \underline{p_2 s_1 s_2} p_1^{-1}$      & $\eqref{P8}$      \\
     & & $s_2 s_1 s_1 s_2$                                                   \\
  \end{longtable}

  \begin{longtable}{cr@{$\ =\ $}lc} \deductionTableWidth
    \eqref{6-R1}: & $r_2 t_2 r_2^{-1}$ 
       & $p_2 \underline{p_1 t_2} p_1^{-1} p_2^{-1}$     & $\eqref{P11}$     \\
     & & $\underline{p_2 t_1} p_2^{-1}$                  & $\eqref{P10}$     \\
     & & $t_1$                                                               \\
  \end{longtable}

  \begin{longtable}{cr@{$\ =\ $}lc} \deductionTableWidth
    \eqref{7-R1}: & $r_2 t_3 r_2^{-1}$ 
       & $p_2 \underline{p_1 t_3} p_1^{-1} p_2^{-1}$     & $\eqref{P10}$     \\
     & & $\underline{p_2 t_3} p_2^{-1}$                  & $\eqref{P11}$     \\
     & & $t_2$ 
  \end{longtable}

  \begin{longtable}{cr@{$\ =\ $}lc} \deductionTableWidth
    \eqref{8-R1}: & $r_2 t_k r_2^{-1}$ 
       & $p_2 \underline{p_1 t_k} p_1^{-1} p_2^{-1}$     & $\eqref{P10}$     \\
     & & $\underline{p_2 t_k} p_2^{-1}$                  & $\eqref{P10}$     \\
     & & $t_k$                                                               \\
  \end{longtable}

  \begin{longtable}{cr@{$\ =\ $}lc} \deductionTableWidth
    \eqref{9-R1}: & $r_2 s_2 r_2^{-1}$
       & $\underline{p_2 p_1 s_2} p_1^{-1} p_2^{-1}$     & $\eqref{P7}$      \\
     & & $s_1$                                                               \\
  \end{longtable}

  \begin{longtable}{cr@{$\ =\ $}lc} \deductionTableWidth
    \eqref{10-R1}: & $r_2 s_k r_2^{-1}$ 
       & $\underline{p_2 p_1 s_k} p_1^{-1} p_2^{-1}$     & $\eqref{P5}^2$    \\
     & & $s_k$
  \end{longtable}

  To deduce \eqref{11-R1} we make use of the following deduction.
  \begin{equation} \label{star} \tag{\ensuremath{\star}}
    \begin{array}{rlc}
      \multicolumn{2}{l}{\underline{p_2 p_1} s_1 s_2 t_3 p_2 p_1} 
                                                                & \eqref{P7} \\
      \qquad = & s_1^{-1} p_2 p_1 \underline{s_2 s_1 s_2} t_3 p_2 p_1 
                                                                & \eqref{P4} \\
             = & s_1^{-1} p_2 p_1 \underline{s_1 s_2 s_1 t_3} p_2 p_1  
                                                  & \eqref{P12}\eqref{P13}^2 \\
             = & s_1^{-1} p_2 p_1 t_1 s_1 \underline{s_2 s_1 p_2} p_1 
                                                                & \eqref{P6} \\
             = & s_1^{-1} p_2 \underline{p_1 t_1 s_1 p_1} s_2 s_1 p_1 
                                                                & \eqref{P9} \\
             = & \underline{s_1^{-1} p_2 s_1} t_1 s_2 s_1 p_1
                                                                & \eqref{P8} \\
             = & s_2 p_1 s_2^{-1} \underline{t_1 s_2} s_1 p_1   & \eqref{P12}\\
             = & s_2 \underline{p_1 t_1 s_1 p_1}                & \eqref{P9}\\
             = & s_2 s_1 t_1                                    & 
    \end{array} 
  \end{equation}

  \begin{longtable}{cr@{$\ =\ $}lc} \deductionTableWidth
    \eqref{11-R1}: & $r_2 s_1 s_2 s_2 s_1 t_1 t_1 r_2^{-1}$
       & $p_2 p_1 s_1 s_2 \underline{s_2 s_1 t_1} t_1 p_1^{-1} p_2^{-1}$
                                                            & $\eqref{P13}^2$\\
     & & $p_2 p_1 s_1 s_2 t_3 \underline{s_2 s_1 t_1 p_1^{-1} p_2^{-1}}$
                                                            & $\eqref{star}$ \\
     & & $\underline{p_2 p_1 s_1 s_2 t_3 p_2 p_1} s_1 s_2 t_3$
                                                            & $\eqref{star}$ \\
     & & $s_2 s_1 \underline{t_1 s_1 s_2} t_3$              & $\eqref{P13}^2$\\
     & & $s_2 s_1 s_1 s_2 t_3 t_3$                                           \\
                                                                             \\
    \eqref{12-R1}: & $r_2 s_3 s_2 s_1 s_1 s_2 s_3 r_2^{-1}$ 
       & $p_2 \underline{p_1 s_3} s_2 s_1 s_1 s_2 s_3 p_1^{-1} p_2^{-1}$
                                                            & $\eqref{P5}$   \\
     & & $p_2 s_3 \underline{p_1 s_2 s_1} s_1 s_2 s_3 p_1^{-1} p_2^{-1}$
                                                            & $\eqref{P6}$   \\
     & & $p_2 s_3 s_2 s_1 \underline{p_2 s_1 s_2} s_3 p_1^{-1} p_2^{-1}$
                                                            & $\eqref{P8}$   \\
     & & $p_2 s_3 s_2 s_1 s_1 s_2 \underline{p_1 s_3} p_1^{-1} p_2^{-1}$
                                                            & $\eqref{P5}$   \\
     & & $\underline{p_2 s_3 s_2} s_1 s_1 s_2 s_3 p_2^{-1}$ & $\eqref{P6}$   \\
     & & $s_3 s_2 \underline{p_3 s_1 s_1} s_2 s_3 p_2^{-1}$ & $\eqref{P5}^2$ \\
     & & $s_3 s_2 s_1 s_1 \underline{p_3 s_2 s_3} p_2^{-1}$ & $\eqref{P8}$   \\
     & & $s_3 s_2 s_1 s_1 s_2 s_3$                                           \\
                                                                             \\
    \eqref{2-R2}: & $r_2 s_1 t_2 s_2 r_2$
       & $p_2 p_1 s_1 \underline{t_2 s_2 p_2} p_1$          & $\eqref{P9}$   \\
     & & $p_2 p_1 s_1 p_2^{-1} s_2 t_2 \underline{p_1}$     & $\eqref{P9}$ \\
     & & $p_2 p_1 s_1 p_2^{-1} s_2 \underline{t_2 s_1^{-1}} t_1^{-1} p_1^{-1} s_1 t_1$
                                                            & $\eqref{P13}$ \\
     & & $p_2 p_1 \underline{s_1 p_2^{-1}} s_2 s_1^{-1} p_1^{-1} s_1 t_1$
                                                            & $\eqref{P6}$   \\
     & & $p_2 p_1 s_2^{-1} p_1^{-1} \underline{s_2 s_1 s_2} s_1^{-1} %
                                          p_1^{-1} s_1 t_1$ & $\eqref{P4}$   \\
     & & $p_2 p_1 s_2^{-1} p_1^{-1} \underline{s_1 s_2 p_1^{-1}} s_1 t_1$
                                                            & $\eqref{P8}$   \\
     & & $p_2 p_1 \underline{s_2^{-1} p_1^{-1} p_2^{-1}} s_1 s_2 s_1 t_1$
                                                            & $\eqref{P7}$   \\
     & & $s_2 s_1 t_1$
  \end{longtable}  

  \begin{longtable}{cr@{$\ =\ $}lc} 
     xxxxxx & xxx & xxxxxxxxxxxxxxxxxxxxxxxxxxxxxxxxx & xxxxxxxxxx \kill
    \eqref{1-R3}: & \multicolumn{3}{l}{$r_1 s_1 s_2 s_3 s_1 s_2 r_1 s_1 
                    s_2 s_3 s_1 s_2 t_2 t_4 r_1 s_2 s_3 s_1 s_2 r_1$}        \\
     & & $p_1 s_1 s_2 s_3 s_1 s_2 p_1 s_1 s_2 s_3 \underline{s_1 s_2
                    t_2} t_4 p_1 s_2 s_3 s_1 s_2 p_1$       & $\eqref{P13}\eqref{P12}$ \\
     & & $p_1 s_1 s_2 s_3 s_1 s_2 p_1 \underline{s_1 s_2 s_3 t_3} s_1
                    s_2 t_4 p_1 s_2 s_3 s_1 s_2 p_1$        & $\eqref{P13}\eqref{P12}^2$ \\
     & & $p_1 s_1 s_2 s_3 \underline{s_1 s_2 p_1 t_4} s_1 s_2 s_3 s_1 
                    s_2 t_4 p_1 s_2 s_3 s_1 s_2 p_1$        & $\eqref{P10}\eqref{P12}^2$ \\
     & & $p_1 \underline{s_1 s_2 s_3 t_4} s_1 s_2 p_1 s_1 s_2 s_3 s_1
                    s_2 t_4 p_1 s_2 s_3 s_1 s_2 p_1$        & $\eqref{P13}^3$ \\
     & & $p_1 t_1 s_1 s_2 s_3 s_1 s_2p_1 s_1 s_2 s_3 s_1 s_2 
                   \underline{t_4 p_1} s_2 s_3 s_1 s_2 p_1$ & $\eqref{P10}$  \\
     & & $p_1 t_1 s_1 s_2 s_3 s_1 s_2 p_1 s_1 \underline{s_2 s_3 s_1
                    s_2 p_1} t_4 s_2 s_3 s_1 s_2 p_1$       & $\eqref{P8}^2$ \\
     & & $p_1 t_1 s_1 s_2 s_3 s_1 s_2 \underline{p_1 s_1 p_3} s_2 s_3
                    s_1 s_2 t_4 s_2 s_3 s_1 s_2 p_1$        & $\eqref{P5}\eqref{P1}$ \\
     & & $p_1 t_1 s_1 s_2 \underline{s_3 s_1} s_2 p_3 p_1 s_1 s_2 s_3
                    s_1 s_2 t_4 s_2 s_3 s_1 s_2 p_1$        & $\eqref{P3}$   \\
     & & $p_1 t_1 s_1 s_2 s_1 \underline{s_3 s_2 p_3} p_1 s_1 s_2 s_3 
                    s_1 s_2 t_4 s_2 s_3 s_1 s_2 p_1$        & $\eqref{P6}$   \\
     & & $p_1 t_1 s_1 \underline{s_2 s_1 p_2} s_3 s_2 p_1 s_1 s_2 s_3 
                    s_1 s_2 t_4 s_2 s_3 s_1 s_2 p_1$        & $\eqref{P6}$   \\
     & & $p_1 t_1 s_1 p_1 s_2 \underline{s_1 s_3} s_2 p_1 s_1 s_2 s_3
                    s_1 s_2 t_4 s_2 s_3 s_1 s_2 p_1$        & $\eqref{P3}$   \\
     & & $\underline{p_1 t_1 s_1 p_1} s_2 s_3 s_1 s_2 p_1 s_1 s_2 s_3
                    s_1 s_2 t_4 s_2 s_3 s_1 s_2 p_1$        & $\eqref{P9}$   \\
     & & $s_1 t_1 s_2 s_3 s_1 s_2 p_1 s_1 s_2 s_3 s_1 s_2 t_4 
                    \underline{s_2 s_3 s_1 s_2 p_1}$        & $\eqref{P8}^2$ \\
     & & $s_1 t_1 s_2 s_3 s_1 s_2 p_1 \underline{s_1 s_2 s_3 s_1 s_2
                    t_4} p_3 s_2 s_3 s_1 s_2$               & $\eqref{P12}^2 \eqref{P13}^3$ \\
     & & $s_1 t_1 s_2 s_3 s_1 s_2 p_1 t_1 s_1 s_2 \underline{s_3 s_1}
                    s_2 p_3 s_2 s_3 s_1 s_2$                & $\eqref{P3}$   \\
     & & $s_1 t_1 s_2 s_3 s_1 s_2 p_1 t_1 s_1 \underline{s_2 s_1 s_3 
                    s_2 p_3} s_2 s_3 s_1 s_2$               & $\eqref{P6}^2$ \\
     & & $s_1 t_1 s_2 s_3 s_1 s_2 \underline{p_1 t_1 s_1 p_1} s_2 s_1
                    s_3 s_2 s_2 s_3 s_1 s_2$                & $\eqref{P9}$   \\
     & & $s_1 t_1 s_2 s_3 s_1 s_2 s_1 \underline{t_1 s_2 s_1} s_3 s_2
                    s_2 s_3 s_1 s_2$                        & $\eqref{P12}\eqref{P13}$ \\
     & & $s_1 t_1 s_2 s_3 s_1 s_2 s_1 s_2 s_1 \underline{t_2 s_3 s_2
                    s_2} s_3 s_1 s_2$                       & $\eqref{P12}\eqref{P13}^2$ \\
     & & $s_1 t_1 s_2 s_3 s_1 s_2 s_1 s_2 s_1 s_3 s_2 s_2 
                    \underline{t_2 s_3 s_1} s_2$            & $\eqref{P12}\eqref{P13}$ \\
     & & $s_1 t_1 s_2 s_3 s_1 s_2 s_1 s_2 s_1 s_3 s_2 s_2 s_3 s_1
                    \underline{t_1 s_2}$                    & $\eqref{P12}$  \\
     & & $s_1 \underline{t_1 s_2 s_3 s_1 s_2} s_1 s_2 s_1 s_3 s_2 s_2 
                    s_3 s_1 s_2 t_1$                        & $\eqref{P12}^2 \eqref{P13}^2$ \\
     & & $s_1 s_2 s_3 s_1 s_2 \underline{t_3 s_1 s_2 s_1} s_3 s_2 s_2 
                    s_3 s_1 s_2 t_1$                        & $\eqref{P12}\eqref{P13}^2$ \\
     & & $s_1 s_2 s_3 s_1 s_2 s_1 s_2 s_1 \underline{t_1 s_3 s_2 s_2 
                    s_3 s_1 s_2} t_1$                       & $\eqref{P12}^4 \eqref{P13}^2$ \\
     & & $s_1 s_2 s_3 s_1 s_2 s_1 s_2 s_1 s_3 s_2 s_2 s_3 s_1 s_2 
                    \underline{t_3 t_1}$                    & $\eqref{P14}$  \\
     & & $s_1 s_2 s_3 s_1 s_2 s_1 s_2 s_1 s_3 s_2 s_2 s_3 s_1 s_2 t_1 t_3$
  \end{longtable}

  \begin{figure}[!htb]
    \centering 
    \def\scalefactor{0.98}  
    $\begin{minipage}[c]{2.5cm} \centering
       \includegraphics[scale=\scalefactor]{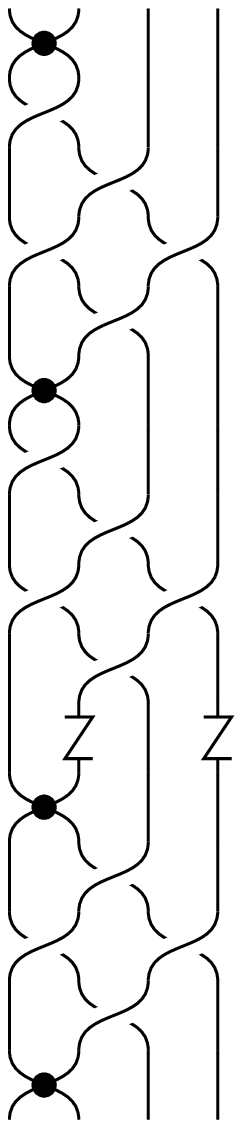}
     \end{minipage}  = \begin{minipage}[c]{2.5cm} \centering
       \includegraphics[scale=\scalefactor]{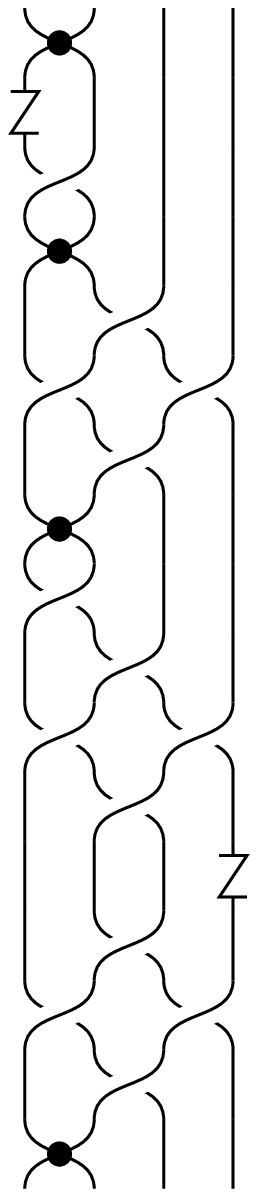}
     \end{minipage} = \begin{minipage}[c]{2.5cm} \centering
       \includegraphics[scale=\scalefactor]{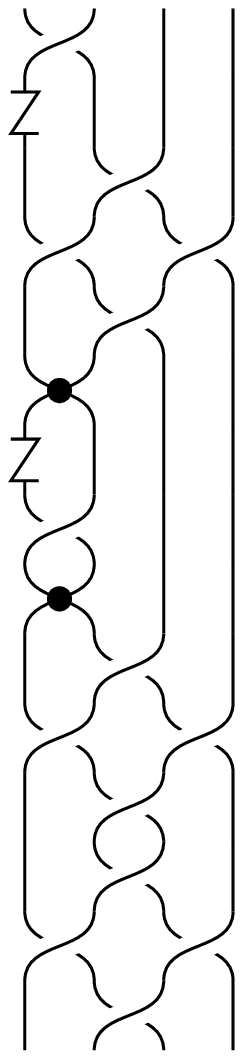}
     \end{minipage} = \begin{minipage}[c]{2.5cm} \centering
       \includegraphics[scale=\scalefactor]{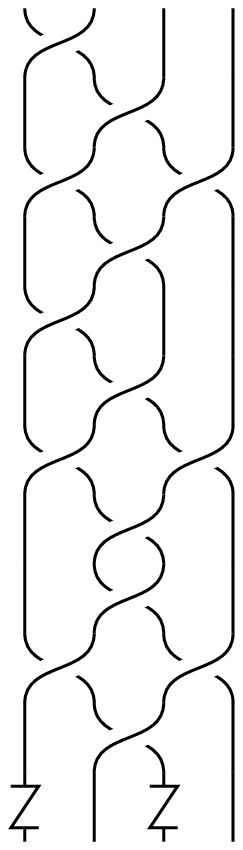}
     \end{minipage}$
    \caption{Deducing the \eqref{1-R3} relation}
    \label{1-R3-pic}
  \end{figure}

  \begin{longtable}{cr@{$\ =\ $}lc} \deductionTableWidth
    \eqref{2-R3}: & $r_1 r_2 s_1 s_2 s_1 t_3 t_2 r_1 r_2$
       & $p_1 p_2 p_1 \underline{s_1 s_2 s_1 t_3} t_2 p_1 p_2 p_1$
                                                & $\eqref{P12}\eqref{P13}^2$ \\
     & & $p_1 p_2 \underline{p_1 t_1 s_1} s_2 s_1 t_2 p_1 p_2 p_1$
                                                           & $\eqref{P9}$    \\
     & & $p_1 p_2 s_1 t_1 \underline{p_1^{-1} s_2 s_1} t_2 p_1 p_2 p_1$
                                                           & $\eqref{P6}$    \\
     & & $p_1 p_2 s_1 t_1 s_2 s_1 p_2^{-1} t_2 \underline{p_1 p_2 p_1}$
                                                           & $\eqref{P7}$    \\
     & & $p_1 p_2 s_1 t_1 s_2 s_1 p_2^{-1} \underline{t_2 p_2} p_1 p_2$
                                                           & $\eqref{P11}$   \\
     & & $p_1 p_2 s_1 \underline{t_1 s_2} s_1 \underline{t_3 p_1} p_2$
                                                  & $\eqref{P12}\eqref{P10}$ \\
     & & $p_1 p_2 s_1 s_2 \underline{t_1 s_1 p_1} t_3 p_2$ & $\eqref{P9}$    \\
     & & $p_1 p_2 \underline{s_1 s_2 p_1^{-1}} s_1 t_1 t_3 p_2$ 
                                                           & $\eqref{P8}$    \\
     & & $p_1 \underline{s_1 s_2 s_1} t_1 t_3 p_2$         & $\eqref{P4}$    \\
     & & $\underline{p_1 s_2 s_1} s_2 t_1 t_3 p_2$         & $\eqref{P6}$    \\
     & & $s_2 s_1 p_2 s_2 \underline{t_1 t_3 p_2}$ 
                                                  & $\eqref{P14}\eqref{P10}$ \\
     & & $s_2 s_1 \underline{p_2 s_2 t_3 p_2} t_1$         & $\eqref{P9}$    \\
     & & $s_2 s_1 s_2 t_2 t_1$                                               \\
                                                                             \\
    \eqref{3-R3}: & $r_2 s_1 t_2 r_1 s_2 s_1 r_1$ 
       & $p_2 \underline{p_1 s_1 t_2 p_1} s_2 s_1 p_1$     & $\eqref{P9}$    \\
     & & $p_2 s_1 \underline{t_1 s_2} s_1 p_1$             & $\eqref{P10}$   \\
     & & $\underline{p_2 s_1 s_2} t_1 s_1 p_1$             & $\eqref{P8}$    \\
     & & $s_1 s_2 \underline{p_1 t_1 s_1 p_1}$             & $\eqref{P9}$    \\
     & & $s_1 s_2 s_1 t_1$
  \end{longtable} 

\end{proof}

\bibliographystyle{plain}
\bibliography{link}

\end{document}